\input amstex
\documentstyle{amsppt}
\NoRunningHeads
\TagsOnRight
\topmatter
\title
On Variational Inequalities with Multivalued
Operators with Semi-Bounded Variation
\endtitle
\author
O. V. Solonoukha\\
National Technical University of Ukraine "KPI", Kiev, Ukraine
\endauthor
\thanks
This work was supported, in part, by the International
Soros Science Education Program (ISSEP) through grant Ø PSU061103.
\endthanks
\endtopmatter

\def\claim#1{\noindent{\bf {#1}.}}
\def\definition#1{\noindent{\bf {#1}.}}
\def\enddefinition{\vskip1pt}
\def\demo#1{\noindent{\it {#1}.}}
\def\enddemo{\qed}
\def\d#1#2{{{\partial#2}\over{\partial x_#1}}}
\def\dnu#1{{{\partial#1}\over{\partial\nu_A}}}
\def\A{\overline{co}A}
\def\[{[\![}
\def\]{]\!]}

\document

In this paper we explore some problems for the steady-state
variational inequalities with multivalued operators (VIMO). As far as we
know, in this variant the term VIMO had been first introduced in
\cite1. The results of this paper with respect to VIMO extend and/or
improve analogous ones from \cite{1-7}. We refused
the regularity conditions of the monotonic  disturbance of the
multivalued mapping (\cite1) and some other
properties of the objects (\cite2). Besides we considered a wider class
of operators with respect to \cite{3,4}. We are  studying
the connections between the class of radially semi-continuous operators
with semi-bounded variation, the class of pseudo-monotone mappings, which
is used earlier (for example, in \cite2) on selector's language, and
the class of monotone mappings.  Moreover, for new class of operators the
property of local boundedness is substituted for a weaker one
with respect to \cite{1,2,5,7,8}.  Also we refuse the condition that $A(y)$
is a convex closet set owing to the forms of support functions.

Let $X$ be a reflexive Banach space, $X^*$ be its topological
dual space, by $\langle\cdot ,\cdot\rangle$ we denote the
dual pairing on $X\times X^*$, $2^{X^*}$ be the totality of all nonempty
subsets of the space $X^*$, $A:  X \to 2^{X^*}$ be a multivalued mapping
with $DomA=\{y\in X:  A(y)\ne\emptyset\}$.
$A:X\to 2^{X^*}$ is called {\it strong} iff $DomA=X$.
Further for simplicity we will consider only strong mappings $A$.
Let us consider the upper and lower
support functions which are associated to $A$:

$\qquad \qquad [A(y),\xi ]_+ = \sup
\limits_{d\in A(y)} \langle d,\xi\rangle,
\qquad  [A(y),\xi ]_- = \inf
\limits_{d\in A(y)} \langle d,\xi\rangle$,

\noindent
and norms:  $\qquad \[Ay\]_+
= \sup\limits_{d\in A(y)}\|d\|_{X^*},\qquad \[Ay\]_- = \inf\limits_{d\in
A(y)}\|d\|_{X^*}$.

We will consider the following VIMO
$$
[A(y),\xi -y]_+  + \varphi (\xi ) -
\varphi (y)\ \geqslant\ \langle f,\xi -y\rangle \ \ \ \ \ \forall \xi\in
\operatorname{dom}\varphi\cap K,\tag1 $$
where $f$ is a fixed element from $X^*$ , $\varphi
:X\to\overline{{\Bbb R}} ={\Bbb R}\bigcup\{+\infty \}$ is a proper convex
lower semi-continuous function, $\operatorname{dom}\varphi
=\{y\in X: \varphi(y)<\infty\}$, $K$ is a convex
weakly closed set.

\definition{Definition}
$L:X\to\overline{{\Bbb R}}$ is called
{\it lower semi-continuous}, if the following is satisfied:
if $X\ni y_n\to y$ in $X$ then
$\varliminf\limits_{n\to\infty} L(y_n)\geqslant L(y)$.

\definition{Definition[1]} Operator $A:X \to 2^{X^*}$ is called

\noindent
a) {\it radially semi-continuous}, if for each $y,\xi,h\in X$ the following
inequality holds:
\centerline{
$\varliminf\limits_{t\to+0}[A(y+t\xi),h]_+ \geqslant [A(y),h]_-$;
}
\noindent
b) {\it operator with semi-bounded variation}, if for each
$R>0$ and arbitrary
$y_1,y_2\in X$ such that $\|y_i\|_X\leqslant R$ ($i=1,2$) the following 
inequality
holds:

\centerline{$ [A(y_1),y_1-y_2]_-\geqslant
[A(y_2),y_1-y_2]_+-C(R;\|y_1-y_2\|'_X), $}
\noindent
where $C:{{\Bbb R}}_+\times{{\Bbb R}}_+\to{{\Bbb R}}$ is continuous,
$\tau^{-1}C(r_1,\tau r_2)\to 0$ as  $\tau\downarrow 0$ for
each $r_1,r_2>0$, $\|\cdot \|'_X$ is a compact norm with respect to the
initial norm $\|\cdot\|_X$;

\noindent
c) {\it coercive operator}, if $\exists y_0\in K$ such that

\centerline{$\|y\|^{-1}_X[A(y),y-y_0]_- \to\infty$ as
$\|y\|_X\to\infty$;}
\noindent
d) {\it locally bounded on $X$}, if for each
$y\in X$ there exist $\varepsilon>0$ and $M>0$ such that
$\[A(\xi)\]_+\leqslant M$ for each $\xi$ such that $\|\xi 
-y\|_X\leqslant\varepsilon$;

\noindent
e) {\it monotone}, if
for each
$y_1,y_2\in X$ the following inequality holds:

\centerline{$ [A(y_1),y_1-y_2]_-\geqslant
[A(y_2),y_1-y_2]_+$.}
\enddefinition

\definition{Definition[2]}  Operator $A:X \to 2^{X^*}$ is called
{\it pseudo-monotone}, if

\noindent
i) the set $A(y)$ is nonempty, bounded, closed and convex at each $y\in X$;

\noindent
ii) $A:F\to 2^{X^*}$ is locally bounded on each
finite-dimensional subspace $F\subset X$;

\noindent
iii) if $y_j\to y$ weakly in $X$, $w_j\in A(y_j)$ and
$\varlimsup\limits_{j\to\infty}\langle w_j,y_j-y\rangle\leqslant 0$, then for
each element $v\in X$ there exists $w(v)\in A(y)$ with the property

\centerline{$
\varliminf\limits_{j\to\infty}\langle
w_j,y_j-v\rangle\geqslant\langle w(v),y-v\rangle$.}
\enddefinition

\definition{Definition[2]}  Operator $A:X \to 2^{X^*}$ is called
{\it generalized pseudo-monotone}, if from $y_j\to y$ weakly in $X$,
$A(y_j)\ni w_j\to w$ $*$-weakly in $X^*$  and
$\varlimsup\limits_{j\to\infty}\langle w_j,y_j-y\rangle\leqslant 0$ it follows
that $w\in A(y)$ and $\langle w_j,y_j\rangle\to\langle w,y\rangle$.
\enddefinition

Each pseudo-monotone operator is generalized pseudo-monotone one (\cite2).

It is easy to see that each monotone operator is an operator with
semi-bounded variation, the next result is connecting the classes of
operators with semi-bounded variation and of pseudo-monotone operators.
Simultaneously, we showed the interconnection between monotone and
pseudo-monotone operators.

\definition{Definition}  Operator $A:X \to 2^{X^*}$ is called
{\it sequentially weakly locally bounded}, if for each
$y\in X$ if $y_n\to y$ weakly in $X$ then there exist a
finite number $N$ and a constant $M>0$ such that $\[A(y_n)\]_+\leqslant M$ for
each $n\geqslant N$.  \enddefinition

\proclaim {Theorem 1}
Let $A$ be a radially semi-continuous operator with semi-bounded variation.
Then $\A$ is pseudo-monotone, locally bounded and sequentially weakly
locally bounded.
\endproclaim

\remark{Remark} It is enough to consider the weaker condition of the
radially semi-conti\-nu\-i\-ty:
$\qquad\qquad\varliminf\limits_{t\to+0}[A(y+t\xi),-\xi]_+ \geqslant
[A(y),-\xi]_-$. \endremark

Let us consider solvability theorems.

\proclaim{Theorem 2}
Let $K$ be a bounded convex weakly closed set and $A:X \to 2^{X^*}$
be a radially semi-continuous operator with semi-bounded
variation. Then for each $f\in X^*$ the solution set of the inequality
$$
[A(y),\xi -y]_+\geqslant\langle f,\xi-y\rangle\ \ \forall\xi\in K\tag2
$$
is nonempty and weakly compact in X. Moreover, there exists the
element $w\in\A(y)$
such that $\qquad\qquad
\langle w,\xi-y\rangle\geqslant\langle f,\xi-y\rangle\quad\forall\xi\in
K$.
\endproclaim
\demo{Proof}
Let us consider the filter ${\Bbb F}$ of the finite-dimensional subspaces
$F$ of $X$. We construct the auxiliary operator $\quad
L_\varepsilon(\lambda,y)=\overline{co}\lbrace(1-\lambda)P_\varepsilon(y)+
\lambda(I_F^*f-I_F^*A(I_Fy))\rbrace,$

\noindent
where $I_F:X\to F$ is the inclusion map, $K_F=K\cap F$,
$P_\varepsilon(y)=[K_F\cap(-N_{K_F}(y))]\setminus B_\varepsilon(y)-y$,
$N_{K_F}(y)$ is the normal cone, $B_\varepsilon(y)=\{\xi\in K_F:
\|\xi-y\|_F<\varepsilon\}$ and $\varepsilon\geqslant 0$ such that
$K_F\setminus B_\varepsilon(y)\ne\emptyset$ for each $y$ from
$\partial K_F$.
We can show that $L_\varepsilon(\lambda,\cdot)$ is upper semi-continuous on
$F$.  By construction for each $y\in\partial K_F$ we have that
$L_\varepsilon(0,y)\cap T_{K_F}(y)\ne\emptyset$,
where $T_{K_F}$ is the tangential cone.
If $\exists y\in\partial K_F$ and $\lambda\in[0,1]$ such that
$0\in L_\varepsilon(\lambda,y)$ then this $y\in\partial K_F$ is a solution
of VIMO on $F$.  Else by Lere--Schauder theorem the inclusion $0\in
L_\varepsilon(0,y)$ has a solution on $intK_F$.  Thus, we have the
bounded sequence $\{y_F\}\subset K$. Using the generalized
pseudo-monotonicity and the sequentially weakly locally boundedness of the
operator $\A$ we can find some limit element $y$ which is a solution of
\thetag2.  \enddemo

\proclaim{Theorem 3}
Let the conditions of Theorem 2 be satisfied without $K$ be a bounded
set. If in this case operator $A$ is coercive, then the statement of
Theorem 2 holds.
\endproclaim
\demo{Proof}
On each bounded set $K_R=K\cap
B_R$ the solution $y_R$ exists, by the coercivity of the operator $\A$
under some $R$ $y_R$ is a solution of \thetag2. \enddemo

From Theorem 3 we can obtain following statement:

\proclaim{Theorem 4}
Let $A:X \to 2^{X^*}$ be a radially semi-continuous coercive operator with
semi-bounded variation. Then $\forall f\in X^*$ the solution set
of the inclusion
\centerline{$\overline{co}A(y)\ni f$}
is nonempty and weakly
compact in $X$.
\endproclaim

Now we consider the based inequality \thetag1
and the corresponding inclusion
$$
\overline{co} A(y) + \partial \varphi (y) \ni f,\tag3
$$
where $\partial\varphi(y)$ is the subdifferential of the function
$\varphi:X\to\overline{\Bbb R}$ at the point $y\in X$.

\proclaim{Proposition}
Each solution of \thetag3 satisfies VIMO
\thetag1. If $y$ is a solution of \thetag1 and belong to
$intK\cap \operatorname{dom}\partial\varphi$, then $y$ is a
solution of \thetag3.
\endproclaim

This simple statement allows to study VIMO \thetag1
using the inclusion \thetag3.

\proclaim{Theorem 5}
Let $A:X \to 2^{X^*}$ be a radially semi-continuous operator with
semi-bounded variation, $\varphi:X\to\overline{\Bbb R}$ be a proper convex
lower semi-conti\-nu\-o\-us function and the following coercivity condition
satisfies:

\centerline{$\exists y_0\in \operatorname{dom}\varphi\cap K$ such that}
\centerline{
$\|y\|_X^{-1}\bigl([A(u,y),y-y_0]_--\varphi(y)\bigr)\to +\infty$
as $\|y\|_X\to\infty.$}
\noindent
Then $\forall f\in X^*$ the
solution set of \thetag1 is nonempty and weakly
compact in $X$.
\endproclaim
\demo{Proof} Let us construct the auxiliary objects:

\centerline{$\widetilde X=X\times \Bbb R,
\quad\widetilde y=(y,\mu)\in \widetilde X,
\quad\widetilde A(\widetilde
y)=(A(y),0)\quad\forall\,\widetilde y\in\widetilde X
$,}
\centerline{$\widetilde K=\{(y,\mu)\in
(K\cap\operatorname{dom}\varphi)
\times\Bbb R\mid \mu\geqslant\varphi(y)\},\quad \widetilde f=(f,-1),$}
\noindent
We can prove that these objects satisfy all conditions of Theorem 2.
Thus, the solution $\widetilde y$ exists and
its first coordinate is a solution of \thetag1.  \enddemo

\claim{Example}
Let us consider the free boundary problem on Sobolev space $W^2_p(\Omega)$,
$p\geqslant 2$:

\centerline{$
-\sum\limits^n_{i,j=1}a_{ij}(x,y,Dy){{\partial^2y}\over {\partial
x_i\partial x_j}}=f$ on $\Omega$,}
\centerline{$y\geqslant 0,\ \  \dnu y\geqslant 0,\ \
y\dnu y =0\quad$  on $\Gamma$,}
\noindent
where $\Omega$ is a sufficiently smooth simple connected
domain of ${\Bbb R}^n$, $\Gamma$ is the bound of $\Omega$, the normal
vector $\nu$ is defined at each $x\in\Gamma$, $\dnu y
=\sum\limits^n_{i,j=1}a_{ij}(x,y,Dy) {{\partial y}\over {\partial
x_j}}cos(x,\nu_i)$, $Dy=(\d 1y,\dots,\d ny)$.  This problem can have not a
classical solution, but we can find a weak solution on
 $W^2_p(\Omega)$.  Let us assume that $a_{ij}(x,y,\xi)$
satisfy the following conditions:  \roster \item for each $y,\xi$ the
functions $a_{ij}$ are continuous with respect to  $x$, \item $\forall
x\in\overline\Omega$ the functions $a_{ij}$ are bounded with respect to
$\xi$ and $y$, and the following estimation holds:  $\quad
|a_{ij}(x,y,\xi_1,\cdots,\xi_n)|\leqslant g(x)+k_0|y|^{p-2}+
\sum\limits^n_{i=1}k_i|\xi_i|^{p-2}$,
where $k_i>0$ ($i=\overline{1,n}$), if $p=2$ then  $g\in C(\Omega)$,
and if $p>2$ then  $g\in L_{q'}(\Omega)$, $q'=p/(p-2)$,
\item  $\sum\limits^n_{i,j=1} a_{ij}(x,y,\xi)\xi_i\xi_j \geqslant
\gamma(R)R$, where $R=|y|+\sum\limits^n_{i=1}|\xi_i|$ and
$\gamma(R)\to+\infty$ as $R\to+\infty$.
\endroster
Then the free boundary problem conforms to the following
inequality:  $$ [A(y),\xi-y]_+=a_1(y,\xi-y)+[A_2(y),\xi-y]_+\geqslant\langle
f,\xi-y\rangle\quad \forall\xi\in K.\tag4 $$ where
$a_1(y,\xi)=\sum\limits^n_{i,j=1}\int\limits_\Omega a_{ij}(x,y,Dy)
{{\partial y}\over {\partial x_j}}{{\partial \xi}\over {\partial x_i}} dx$,
$A_2(y)=\bigl({{\partial }\over {\partial x_i}} a_{ij}(x,y,Dy)\bigr)
{{\partial \xi}\over {\partial x_j}}$,
${{\partial a_{ij}}\over {\partial x_i}}$ is the subdifferential
of $a_{ij}$, $K=\{y\in W^2_p(\Omega):y_{|\Gamma}\geqslant 0\}$ is a convex
weakly closed set.  We can prove that $A$ is a radially semi-continuous
coercive operator with semi-bounded variation. Thus, the inequality
\thetag4 has a solution.

\bigskip
{\smc
O. V. Solonoukha, postgraduate student, National Technical
University of Uk\-ra\-i\-ne "Kiev Polytechnic Institute", Chair of
Mathematical Modelling of Economical Systems, pr. Pobedy, 37, Kiev, Ukraine.
}

\enddocument